\documentclass[10pt,a4paper]{article}
\usepackage{amsmath,amssymb,amsthm}
\usepackage[utf8]{inputenc}

\usepackage[usenames]{color}

\usepackage{color}

\definecolor{c20}{rgb}{0.,0.7,0.}
\definecolor{c30}{rgb}{0.,0.,1.}
\definecolor{c40}{rgb}{1,0.1,0.7}
\definecolor{c50}{rgb}{1,0,0}

\def\pE#1{\textcolor{black}{#1}}

\def\cE#1{\textcolor{black}{#1}}
\def\cE0#1{\textcolor{black}{#1}}
\def\cK#1{\textcolor{black}{#1}}
\def\cN#1{\textcolor{black}{#1}}
\newcommand{\tb}[1]{{\textcolor{black}{#1}}}
\newcommand{\tr}[1]{{\textcolor{black}{#1}}}
\newcommand{\tp}[1]{{\textcolor{black}{#1}}}
\newcommand{\ve}{\varepsilon}

\newcommand{\BQN}{\begin{eqnarray}}
\newcommand{\EQN}{\end{eqnarray}}
\newcommand{\BQNY}{\begin{eqnarray*}}
\newcommand{\EQNY}{\end{eqnarray*}}

\newcommand{\rr}{\mathbb R}

\newtheorem{tm}{Theorem}
\newtheorem{Prop}{Proposition}
\newtheorem{lm}{Lemma}

\newtheorem{uw}{Remark}

\newenvironment{pf}{\textsc{Proof.}}{\hfill$\blacksquare$\\\par}

\begin{document}


\title{Extremes of homogeneous Gaussian random fields}
\author{{\small {Krzysztof D\c{e}bicki$^{1}$}, Enkelejd Hashorva$^{2}$, Natalia Soja-Kukie\l a$^{3}$
}\\
{\small\it 1. Mathematical Institute, University of Wroc\l aw, pl. Grunwaldzki 2/4, 50-384 Wroc\l aw, Poland,}\\
{\small\it 2. Department of Actuarial Science, Faculty of Business and Economics,}\\
{\small\it  University of Lausanne, UNIL-Dorigny 1015 Lausanne, Switzerland}\\
{\small\it 3. Nicolaus Copernicus University, ul. Chopina 12/18,
87-100 Toru\'n, Poland}\\
}
\maketitle
\abstract
Let $\{X(s,t):s,t\geqslant 0\}$ be a centered homogeneous Gaussian field with a.s. continuous sample paths
and \cE0{correlation} function $r(s,t)=Cov(X(s,t),X(0,0))$ such that
\[r(s,t)=1-|s|^{\alpha_1}-|t|^{\alpha_2}+o(|s|^{\alpha_1}+|t|^{\alpha_2}), \quad \cE0{s,t \to 0},\]
with $\alpha_1,\alpha_2\in(0,2],$ \cE0{and}
$r(s,t)<1$ for $(s,t)\neq(0,0)$. \cE0{In this contribution we derive an exact asymptotic expansion (as $u\to \infty$) of
$$\mathbb{P}\left(\sup_{(s n_1(u),t n_2(u))\in\left[0,x\right]\times\left[0,y\right]}X(s,t)\leqslant u\right),$$
%
where
$n_1(u)n_2(u)=u^{2/\alpha_1+2/\alpha_2}\Psi(u)$,
which holds uniformly for $(x,y) \in [ \pE{A} ,  \pE{B} ]^2$ with $ \pE{A} , \pE{B} $ two positive constants
 and $\Psi$ the survival function of an $N(0,1)$ random variable.
\cK{We apply our findings to the analysis of asymptotics of
extremes of homogeneous Gaussian fields over
\tb{more complex parameter sets and
a ball of random radius. Additionally}
we determine the extremal index of
the discretised random field determined by $X(s,t)$.}}

\pE{Key words: Gaussian random fields; supremum; tail asymptoticy; extremal index; Berman condition; strong dependence.}

\section{Introduction}
One of the seminal results in extreme value theory of Gaussian processes is the asymptotic behaviour of the distribution of supremum of a centered stationary Gaussian process $\{\cE0{X}(t):t\geqslant 0\}$ with \cE0{correlation function} satisfying
\begin{equation}\label{A1}
\cE0{r(t)}=Cov(\cE0{X}(t),\cE0{X}(0))=1-|t|^\alpha+o(|t|^\alpha) \text{ as } t\to 0 \text{ with } \alpha\in(0,2],
\end{equation}
over intervals of length proportional to
$$\mu(u)=P\left(\sup_{t\in[0,1]}X(t)>u\right)^{-1}(1+o(1)),$$
see, e.g., Leadbetter et al. \cite[Theorem 12.3.4]{LEAD}, Arendarczyk and D\c{e}bicki \cite[Lemma 4.3]{AREN}, Tan and Hashorva \cite[Lemma 3.3]{TAN}.
The following theorem gives a preliminary result concerning the aforementioned asymptotics.
\begin{tm}\label{TW1WYM}
Let $\{\cE0{X}(t):\;t\geqslant 0\}$ be a centered stationary Gaussian process that satisfies (\ref{A1}), and let $0<A_0<A_\infty<\infty$ and $x>0$ be arbitrary constants.
\tb{
If $\;\cE0{r}(t)\log t\to r \in [0,\infty)$ as $t\to\infty$, then
$$P\left(\sup_{t\in[0,x\mu(u)]}X(t)\leqslant u\right)\to E\left(\exp\left(-x\exp(-r+\sqrt{2r}\mathcal{W})\right)\right)\in (0,\infty),$$
as $u\to\infty$, uniformly for $x\in[A_0,A_\infty]$, with  $\mathcal{W}$ an $N(0,1)$ random variable.
}
\end{tm}

The main goal of this paper is to \cE0{derive} an analogue of the above \tr{result} for Gaussian random fields;
see part (i) of Theorem \ref{LEM2WYM} which constitutes a $2$-dimensional counterpart of Theorem \ref{TW1WYM}.

\cK{As an application of \cE0{our findings}, in Section
\ref{RESULTS} we investigate asymptotics of the tail of supremum
of a homogeneous Gaussian field over a
\tb{parameter sets that are approximable by simple sets (part (ii) of Theorem \ref{LEM2WYM}) and}
a ball of random radius.
Additionally
we analyze the existence of the {\it
extremal index} for discrete-parameter fields associated with
homogeneous Gaussian fields with covariance structure satisfying
some regularity conditions; see Proposition \ref{prop.1}.}


\section{Preliminaries}\label{PRE}
Let $\{X(s,t):s,t\geqslant 0\}$ be a centered homogeneous Gaussian field with a.s. continuous sample paths and \cE0{correlation} function $r(s,t)=Cov(X(s,t),X(0,0))$ such that
\begin{itemize}
\item[\textbf{A1}:]
$r(s,t)=1-|s|^{\alpha_1}-|t|^{\alpha_2}+o(|s|^{\alpha_1}+|t|^{\alpha_2})$ as $s,t\to 0$ with $\alpha_1,\alpha_2\in(0,2]$;
\item[\textbf{A2}:]
$r(s,t)<1$ for $(s,t)\neq(0,0)$;
\item[\textbf{A3}:]
$\sup_{(s,t)\in \mathcal{S}(0,d)}|r(s,t)\log d - r|\to 0$ as $d\to\infty$, with $r\in [0,\infty)$,
\end{itemize}
where $\mathcal{S}(0,d)$ denotes the sphere of center $(0,0)$ and radius $d>0$ in $\mathbb{R}^2$ with Euclidean metric.

We distinguish two separate families of \cE0{Gaussian fields}
\begin{itemize}
\item {\it weakly dependent fields}, satisfying \textbf{A3} with $r=0$,
\item {\it strongly dependent fields}, satisfying \textbf{A3} with $r\in(0,\infty)$.
\end{itemize}

Let $\mathcal{H}_\alpha$ denote the Pickands constant (see \cite{PICKANDS}), i.e.,
$$\mathcal{H}_\alpha:=
\lim_{T\to\infty}\frac{E\exp\left(\max_{0\leqslant t\leqslant T}\chi(t) \right)}{T}$$
where $\chi(t)=B_{\alpha/2}(t)-|t|^{\alpha}$, with $\{B_{\alpha/2}(t):t\geqslant 0\}$
being a fractional Brownian motion with Hurst parameter $\alpha/2\in(0,1]$. \cE0{We note in passing that
$\mathcal{H}_\alpha$ appears for the first time in Pickands theorem \tr{\cite{PICKANDS}};
a correct proof of that theorem is first given in Piterbarg \cite{Pit72}.}

For a standard normal random variable $\mathcal{W}$ we write $\Phi(u)=P(\mathcal{W}\leqslant u)$,
$\Psi(u)=P(\mathcal{W}>u)$.
Recall that
$$\Psi(u)=\frac{1}{\sqrt{2\pi}u}\exp(-u^2/2)(1+o(1)),
\quad\text{as }u\to\infty.$$

Following Piterbarg \cite[Theorem 7.1]{PITERBARG} we recall that for a centered stationary Gaussian field $\{X(s,t)\}$ satisfying \textbf{A1}, \textbf{A2}, for arbitrary $g,h\in(0,\infty)$,
\begin{equation}\label{ASYMPTOTICS}
P\left(\max_{(s,t)\in[0,g]\times[0,h]}X(s,t)>u\right)=\mathcal{H}_{\alpha_1}\mathcal{H}_{\alpha_2} gh u^{2/\alpha_1}u^{2/\alpha_2}\Psi(u)(1+o(1)),
\end{equation}
as $u\to\infty$.

\tp{Let $m_1(u)\to\infty$ and $m_2(u)\to\infty$ be functions such that
$$m_1(u)=a_1(u)/\sqrt{\Psi(u)}\quad\text{and}\quad m_2(u)=a_2(u)/\sqrt{\Psi(u)}$$ for some positive \pE{functions} $a_1(u),a_2(u)$ satisfying $a_1(u)a_2(u)=(\mathcal{H}_{\alpha_1}\mathcal{H}_{\alpha_2}u^{2/\alpha_1}u^{2/\alpha_2})^{-1}$, $\log a_1(u) =o(u^2)$ and $\log a_2(u) =o(u^2)$.
We note that then
$$m(u):=m_1(u)m_2(u)=P\left(\max_{(s,t)\in[0,1]^2}X(s,t)>u\right)^{-1}(1+o(1)),$$
as $u\to\infty$.}

\cK{
By $\mathcal{B}(0,x)$ we denote a ball in $\rr^2$ of center at $(0,0)$ and radius $x$.
}

\section{Main results}\label{RESULTS}
The aim of this section is to prove the following $2$-dimensional counterpart of Theorem \ref{TW1WYM}. Recall that $\mathcal{W}$ denotes an $N(0,1)$ random variable.
For a given Jordan-measurable set \pE{$\mathcal{E}\subset\rr^2$ with Lebesgue measure
${\rm mes}(\mathcal{E})>0$} let
\tb{$\mathcal{E}_u:=\{(x,y):(x/ m_1(u),y /m_2(u))\in\mathcal{E}\}$.} \pE{One interesting example is $\mathcal{E}_u=
\left[0,xm_1(u)\right]\times\left[0,ym_2(u)\right]$ for $x,y$ positive, hence $\mathcal{E}=[0,x]\times [0,y]$ and
${\rm mes}(\mathcal{E})=xy$. For such $\mathcal{E}_u$ we shall show below an approximation
which  holds uniformly on compact intervals of $(0,\infty)^2$.  If the structure of the set is not specified, considering thus the supremum of a Gaussian field over
some general measurable set $\mathcal{T}_u\subset \rr^2$ an $\epsilon$-net $(\mathcal{L}_\ve,\mathcal{U}_\ve)$ approximation of $\mathcal{T}_u$ will be assumed. Specifically, the
$\epsilon$-net $(\mathcal{L}_\ve,\mathcal{U}_\ve)$ here means that for any $\ve>0$ there exist two sets $\mathcal{L}_\ve $ and $\mathcal{U}_\ve$  which are
{\it simple sets}  (i.e., finite sums of disjoint
rectangles of the form $[a_1,b_1)\times[a_2,b_2)$) such that
\BQN\label{epsA}
\lim_{\ve \downarrow 0} {\rm mes}(\mathcal{L}_\ve)=\lim_{\ve \downarrow 0}
{\rm mes}(\mathcal{U}_\ve)=c\in(0,\infty)
\EQN
and
\BQNY
\mathcal{L}_{\varepsilon,u}=\{(x,y):(x/m_1(u),y/m_2(u))\in\mathcal{L_\varepsilon}\} \subset \mathcal{T}_u \subset
\mathcal{U}_{\varepsilon,u}=\{(x,y):(x/m_1(u),y/m_2(u))\in\mathcal{U_\varepsilon}\}\subset \rr^2.
\EQNY
Next we formulate our main results for these two cases.}

\begin{tm}\label{LEM2WYM}
Let $\{X(s,t):\;s,t\geqslant 0\}$ be a centered homogeneous Gaussian field with covariance function that satisfies \textbf{A1}, \textbf{A2} and \textbf{A3} with $r\in[0,\infty)$.
\tb{ Then,\\
(i) for each $0< \pE{A} < \pE{B} <\infty$,
$$\mathbb{P}\left(\sup_{(s,t)\in\left[0,xm_1(u)\right]\times\left[0,ym_2(u)\right]}X(s,t)\leqslant u\right)
\to
\mathbb{E}\big(\exp(-xy\exp(-2r+2\sqrt{r}\mathcal{W}))\big),$$
as $u\to\infty$, uniformly for $(x,y)\in[ \pE{A} , \pE{B} ]^2$.\\
(ii) 
for $\mathcal{T}_u\subset \rr^2,u> 0$ such that there exists an $\epsilon$-net $(\mathcal{L}_\ve,\mathcal{U}_\ve)$ satisfying \eqref{epsA}
\[
\mathbb{P}\left(\sup_{(s,t)\in \mathcal{T}_u}X(s,t)\leqslant u\right)
\to
\mathbb{E}\big(\exp(- c\exp(-2r+2\sqrt{r}\mathcal{W}))\big),\ {\rm as} \ u\to\infty.
\]
}
\end{tm}
The complete proof of Theorem \ref{LEM2WYM} is given in Section \ref{s.proof.th2}.

\begin{uw}\label{r.2}
\tb{
Following the same reasoning as given in the proof of Theorem \ref{LEM2WYM}, assuming that
{\bf A1-A3} holds,
for each $0< \pE{A} < \pE{B} <\infty$, we have}
\begin{eqnarray}
\mathbb{P}\left(\sup_{(s,t)\in \mathcal{B}\left(0,x\sqrt{m(u)}\right) }X(s,t)\leqslant u\right)
\to
\mathbb{E}\big(\exp(-\pi x^2\exp(-2r+2\sqrt{r}\mathcal{W}))\big),
\end{eqnarray}
as $u\to\infty$, uniformly for $x\in[ \pE{A} , \pE{B} ]$;
$\mathcal{B}(0,x)$ is a ball in $\rr^2$ of center at $(0,0)$ and radius $x$.
\end{uw}

\section{Applications} \label{s. examples}
In this section we apply results of Section \ref{RESULTS} to the analysis
of the asymptotic properties of supremum of a Gaussian field over a random parameter set and
to the analysis of dependance structure of homogeneous Gaussian fields.

\subsection{Extremes of homogeneous Gaussian fields over a random parameter set}
In this section we analyze asymptotic properties of the tail distribution of
$\sup_{(s,t)\in \mathcal{B}(0,T)}X(s,t)>u)$, where $T$ is a nonnegative, independent of
$X$ random variable.
One-dimensional counterpart of this problem was recently analyzed in
\cite{AREN} and \cite{TAN}.

\begin{Prop}\label{p.random}
Let $\{X(s,t):\;s,t\geqslant 0\}$ be a centered homogeneous Gaussian field with covariance
function that satisfies \textbf{A1-A3} with $r\in[0,\infty)$, and
let $T$ be an independent of $X$ nonnegative random variable.\\
(i) If $ET^2<\infty$, then, as $u\to\infty$,
\[P\left(\sup_{(s,t)\in \mathcal{B}(0,T)}X(s,t)>u\right)=\tp{\pi ET^2 \mathcal{H}_{\alpha_1}\mathcal{H}_{\alpha_2} u^{2/\alpha_1}u^{2/\alpha_2}\Psi(u)(1+o(1))}.\]
(ii) If $T$ has a regularly varying survival function at infinity with index $\lambda<2$, then as $u\to\infty$,
\[P\left(\sup_{(s,t)\in \mathcal{B}(0,T)}X(s,t)>u\right)=2\pi \mathcal{C} P(T>\sqrt{m(u)})(1+o(1)),
\]
where $\mathcal{C}=
\int_{0}^{\infty}
x^{1-\lambda}E\big(\exp(-\pi x^2\exp(\mathcal{V}_r)+\mathcal{V}_r)\big)dx
$ and $\mathcal{V}_r=2\sqrt{r}\mathcal{W}-2r$.\\
(iii) If $T$ is slowly varying at $\infty$, then, as $u\to\infty$,
\[P\left(\sup_{(s,t)\in \mathcal{B}(0,T)}X(s,t)>u\right)=P(T>\sqrt{m(u)})(1+o(1)).\]
\end{Prop}
The proof of Proposition \ref{p.random} is given in Section \ref{sub3}.

\subsection{Extremal indices for homogeneous Gaussian fields}
Following \tr{\cite{SOJAK}}, we say that \tr{$\theta\in (0,1]$} is the {\it extremal index}
of a homogeneous
discrete-parameter stationary random field
$\{X_{j,k}:j,k=1,2,\ldots\}$, if
\begin{equation}\label{DISDEF}
P\left(\max_{j\leqslant a_n,\;k\leqslant b_n}X_{j,k} \leqslant z_n\right) - P(X_{1,1}\leqslant z_n)^{a_nb_n\cdot\theta} \to 0,
\end{equation}
as $n\to\infty$, for each sequence $(z_n)\subset\mathbb{R}$
and all sequences $(a_n),(b_n)\subset \mathbb{N}$
such that $a_n\to\infty$ and $b_n\to\infty$, as $n\to\infty$, and \tr{$1/C\leqslant a_n/b_n \leqslant C$} for some constant \tr{$C>0$}.
The notion of extremal index $\theta$ originated in investigations concerning
relationship between the dependence structure of discrete-parameter stationary sequences of random variables
and their extremal behaviour \cite{LEAD2, LEAD};
see also \cite{OBRIEN,HSING,HSING2,JAKUBOW,LEAD3,WEISS}.

For a given centered homogeneous Gaussian field $\{X(s,t):s,t\ge0\}$ that satisfies
{\bf A1-A3} introduce
a discrete-parameter random field
$\{\widetilde{X}_{j,k}: j,k=1,2,\ldots\}$, with
\[\widetilde{X}_{j,k}:=\sup_{(s,t)\in [j-1,j]\times [k-1,k]} X(s,t).\]
The following proposition  points \cE0{out} how the difference in the dependance structure
between weakly- and strongly-dependant Gaussian fields influences the existence of
the extremal index \tr{of the associated field $\{\widetilde{X}_{j,k}\}$.}

\begin{Prop}\label{prop.1}
Assume that {\bf A1-A3} holds for a centered
homogeneous Gaussian field $\{X(s,t):s,t\geqslant 0\}$. \\
(i) If $r=0$, then the extremal index of
$\{\widetilde{X}_{j,k}: j,k=1,2,\ldots\}$ equals to 1.\\
(ii) If $r>0$, then
$\{\widetilde{X}_{j,k}: j,k=1,2,\ldots\}$ does not have an extremal index.
\end{Prop}
The proof of Proposition \ref{prop.1} is deferred to Section \ref{SUB2}.


\section{Proofs}\label{s.proofs}

Before we prove Theorem \ref{LEM2WYM}, we need some auxiliary results.
The first one is a $2$-dimensional version of Lemma 12.2.11 in \cite{LEAD}.

\begin{lm}\label{LEMAT1}
Assume that \textbf{A1}, \textbf{A2} hold and
$q_1=q_1(u)=au^{-2/\alpha_1}$, $q_2=q_2(u)=au^{-2/\alpha_2}$
for some $a>0$. Then for any $x,y\geqslant 0$, $g,h>0$ and rectangle $I=(x,y)+[0,g]\times[0,h]$, as $u\to\infty$,
$$P\left(X(jq_1,kq_2)\leqslant u;\; (jq_1,kq_2)\in I\right)-P\left(X(s,t)\leqslant u;\; (s,t)\in I\right)\leqslant \frac{gh\rho(a)}{m(u)}+o\left(\frac{1}{m(u)}\right),$$
where $\rho(a)\to 0$ as $a\to 0$.
\end{lm}
\begin{pf}
From the homogeneity of the field $\{X(s,t)\}$ we conclude that
\begin{eqnarray*}
\lefteqn{0\leqslant P\left(X(jq_1,kq_2)\leqslant u;\; (jq_1,kq_2)\in I\right)-P\left(X(s,t)\leqslant u;\; (s,t)\in I\right)}\\
&& \leqslant ([g/q_1]+[h/q_2]+1)P(X(0,0)>u)+ P\left(X(jq_1,kq_2)\leqslant u;\; (jq_1,kq_2)\in [0,g]\times[0,h] \right)\\
&& -P\left(X(s,t)\leqslant u;\; (s,t)\in [0,g]\times[0,h] \right).
\end{eqnarray*}
Then there exists a constant $K$ such that
$$
([g/q_1]+[h/q_2]+1)P(X(0,0)>u)m(u) \leqslant\frac{K(u^{2/\alpha_1}+u^{2/\alpha_2}) \Psi(u)}{\mathcal{H}_{\alpha_1}\mathcal{H}_{\alpha_2} u^{2/\alpha_1}u^{2/\alpha_2}\Psi(u)}
,
$$
which implies that $([g/q_1]+[h/q_2]+1)P(X(0,0)>u)=o\left(\frac{1}{m(u)}\right)$, as $u\to\infty$.

Let $T>0$ be given. We divide the set $[0,g]\times[0,h]$ into small rectangles with the side-lengths $q_1T$ and $q_2T$ in the following way
\begin{eqnarray*}
\Delta_{1,1}&:=&[0,q_1T]\times[0,q_2T],\\
\Delta_{l,m}&:=&((l-1)q_1T,(m-1)q_2T)+\Delta_{1,1},
\end{eqnarray*}
for $l =1,\ldots, \left\lfloor\frac{g}{q_1T}\right\rfloor$ and $m=1,\ldots,\left\lfloor\frac{h}{q_2T}\right\rfloor$.
Then we have that
\begin{eqnarray}\label{SZACOWANIE}
\lefteqn{P(X(jq_1,kq_2)\leqslant u;\; (jq_1,kq_2)\in[0,g]\times[0,h])-P(X(s,t)\leqslant u;\; (s,t)\in[0,g]\times[0,h])}\nonumber\\
&& \leqslant P\left(\sup_{(s,t)\in[0,g]\times[0,h]} X(s,t)> u\right) -
\sum_{l=1}^{\left\lfloor\frac{g}{q_1T}\right\rfloor}\sum_{m=1}^{\left\lfloor\frac{h}{q_2T}\right\rfloor}
 P\left(\max_{(jq_1,kq_2)\in \Delta_{l,m}} X(jq_1,kq_2)>u\right)\nonumber\\
&& +\sum_{(l,m)\neq(l',m')}P\left( \max_{(jq_1,kq_2)\in \Delta_{l,m}} X(jq_1,kq_2)>u, \max_{(jq_1,kq_2)\in \Delta_{l',m'}} X(jq_1,kq_2)>u\right).
\end{eqnarray}
From \cite[Lemma 7.1]{PITERBARG}, as $u\to\infty$,
\begin{equation}\label{1}
P\left(\sup_{(s,t)\in[0,g]\times[0,h]}X(s,t)>u\right) =\mathcal{H}_{\alpha_1}\mathcal{H}_{\alpha_2}ghu^{2/\alpha_1}u^{2/\alpha_2}\Psi(u)(1+o(1)).
\end{equation}
Moreover, by homogeneity of $X(\cdot,\cdot)$,
\begin{eqnarray}
\sum_{l=1}^{\left\lfloor\frac{g}{q_1T}\right\rfloor}\sum_{m=1}^{\left\lfloor\frac{h}{q_2T}\right\rfloor} P\left(\max_{(jq_1,kq_2)\in \Delta_{l,m}} X(jq_1,kq_2)>u\right)\sim \frac{gh u^{2/\alpha_1}u^{2/\alpha_2}}{a^2T^2}P\left(\max_{(jq_1,kq_2)\in \Delta_{1,1}} X(jq_1,kq_2)>u\right).
\label{dub.1}
\end{eqnarray}
We focus on the asymptotics of $P\left(\max_{(jq_1,kq_2)\in \Delta_{1,1}} X(jq_1,kq_2)>u\right)$.
\tb{
Following
line-by-line the idea of the proof of Lemma D.1 in \cite{PITERBARG} we have
\begin{eqnarray*}
\lefteqn{P\left(\max_{(jq_1,kq_2)\in \Delta_{1,1}} X(jq_1,kq_2)>u\right)}\\
&&\sim
\Psi(u)\int_{-\infty}^{\infty}e^{w-w^2/(2u^2)}P\left(\max_{(ja,ka)\in[0,aT]^2} \chi_u(ja,ka)>w) \;\Big|\; X(0,0)=u-\frac{w}{u}\right)dw,
\\
&&\sim
\Psi(u)
H_{\alpha_1}(T,a) H_{\alpha_2}(T,a),
\end{eqnarray*}
where
$H_{\alpha_i}(T,a):=E\exp\left(\max_{j\in [0,T]} B_{\alpha_i/2}(ja)-|ja|^{\alpha_i}\right),$
with
$B_{\alpha_i/2}(\cdot)$
being a fractional Brownian motion with Hurst parameter $\alpha_i/2$
for $i=1,2$ (see also (12.2.6) in proof of \cite[Lemma 12.2.11]{LEAD}).
}

The above implies that, by (\ref{dub.1}),
\begin{eqnarray}
\lefteqn{
\sum_{l=1}^{\left\lfloor\frac{g}{q_1T}\right\rfloor}\sum_{m=1}^{\left\lfloor\frac{h}{q_2T}\right\rfloor} P\left(\max_{(jq_1,kq_2)\in \Delta_{l,m}} X(jq_1,kq_2)>u\right)
}\nonumber\\
&=& ghu^{2/\alpha_1}u^{2/\alpha_2}\Psi(u)
\left(\frac{H_{\alpha_1}(T,a)}{aT}\right)
\left(\frac{H_{\alpha_2}(T,a)}{aT}\right)(1+o(1))\label{2}
\end{eqnarray}
as $u\to\infty$.

In the next step we prove that the double sum that appears in (\ref{SZACOWANIE}) is negligible, i.e., it is $o\left(\frac{1}{m(u)}\right)$.
Indeed, notice that
\begin{eqnarray}\label{3}
\lefteqn{\sum_{(m,l)\neq(m',l')}P\left( \max_{(jq_1,kq_2)\in \Delta_{m,l}} X(jq_1,kq_2)>u, \max_{(jq_1,kq_2)\in \Delta_{m',l'}} X(jq_1,kq_2)>u\right)}\nonumber\\
&& \leqslant \sum_{(m,l)\neq(m',l')}P\left( \sup_{(s,t)\in \Delta_{m,l}} X(s,t)>u, \sup_{(s,t)\in \Delta_{m',l'}} X(s,t)>u\right)=o\left(\frac{1}{m(u)}\right),
\end{eqnarray}
where (\ref{3}) follows from the proof of \cite[Lemma 6.1]{PITERBARG}.

Now, combining (\ref{1}), (\ref{2}) and (\ref{3}), we conclude that for any $T>0$ and $a>0$ it holds that
\begin{eqnarray*}
\lefteqn{P\left(X(jq_1,kq_2)\leqslant u;\; (jq_1,kq_2)\in[0,g]\times [0,h]\right)-P\left(X(s,t)\leqslant u;\; (s,t)\in[0,g]\times[0,h]\right)}\\
&&\leqslant gh u^{2/\alpha_1}u^{2/\alpha_2}\Psi(u)
\left(\mathcal{H}_{\alpha_1}\mathcal{H}_{\alpha_2}
-\left(\frac{H_{\alpha_1}(T,a)}{aT}\right)\cdot\left(\frac{H_{\alpha_2}(T,a)}{aT}\right)\right)(1+o(1))\\
&&=
gh\frac{1-\left(\frac{H_{\alpha_1}(T,a)}{aT}\cdot\frac{H_{\alpha_2}(T,a)}{aT}\right) \mathcal{H}_{\alpha_1}^{-1}\mathcal{H}_{\alpha_2}^{-1}}{m(u)}
+o\left(\frac{1}{m(u)}\right).
\end{eqnarray*}
Finally, using that
\begin{equation}
\lim_{a\to 0}\lim_{T\to\infty} \frac{H_\alpha(T,a)}{aT}=\mathcal{H}_\alpha,\nonumber
\end{equation}
see e.g. \cite[Lemmas 12.2.4(i),12.2.7(i)]{LEAD},
the thesis of the lemma is satisfied with
$$\rho(a):=1-\lim_{T\to\infty}\left(\frac{H_{\alpha_1}(T,a)}{aT} \cdot\frac{H_{\alpha_1}(T,a)}{aT}\right) \mathcal{H}_{\alpha_1}^{-1}\mathcal{H}_{\alpha_2}^{-1}.$$
This completes the proof.
\end{pf}

Let
\begin{eqnarray}
\rho_T(s,t) &:=& \left\{ \begin{array}{ll}
1,
& 0\leqslant \cN{\max(|s|,|t|)} < 1;\\
|r(s,t)-\frac{r}{\log T}|,
& 1\leqslant \cN{\max(|s|,|t|)} \leqslant T,\end{array} \right. \label{FUN}
\\
\varrho_T(s,t) &:=& \left\{ \begin{array}{ll}
|r(s,t)|+(1-r(s,t))\frac{r}{\log T},
& 0\leqslant \cN{\max(|s|,|t|)} < 1;\\
\frac{r}{\log T},
& 1\leqslant \cN{\max(|s|,|t|)} \leqslant T.\end{array} \right. \label{FUN1}
\end{eqnarray}
The next lemma combines a $2$-dimensional counterpart of Lemma 12.3.1 in \cite{LEAD},
for weakly dependent fields, and Lemma 3.1 in \cite{TAN} for strongly dependent fields.
\begin{lm}\label{LEMACIK2}
Let $\varepsilon>0$ be given. Let $q_1=q_1(u)=au^{-2/{\alpha_1}}$ and $q_2=q_2(u)=au^{-2/{\alpha_2}}$. Suppose that $T_1=T_1(u)\sim \tau m_1(u)$ and $T_2=T_2(u)\sim \tau m_2(u)$ for some $\tau>0$, as $u\to\infty$.
Then, providing that conditions \textbf{A1}, \textbf{A2} and \textbf{A3} with $r\in[0,\infty)$ are fulfilled,
$$
\frac{T_1T_2}{q_1q_2} \!\!
\sum_{(jq_1,kq_2)\in \cN{[-T_1,T_1]\times[-T_2,T_2]-(-\varepsilon,\varepsilon)^2}}
\rho_{T_{\max}} \! (jq_1,kq_2)\exp\left(-\frac{u^2}{1+\max\big(|r(jq_1,kq_2)|,\varrho_{T_{\max}} \! (jq_1,kq_2)\big)}\right)\to 0,
$$
as $u\to\infty$, where $T_{\max}=\max(T_1,T_2)$.
\end{lm}

\begin{pf}
Let $T_1(u)\sim \tau m_1(u)$ and $T_2(u)\sim \tau m_2(u)$ for some $\tau>0$, as $u\to\infty$.
Then,
$$\log (T_1 T_2) +\log\left(\frac{\mathcal{H}_{\alpha_1}\mathcal{H}_{\alpha_2}}{\sqrt{2\pi}}\right) +\left(\frac{2}{\alpha_1}+\frac{2}{\alpha_2}-1\right)\log u -\frac{u^2}{2}
\to 2\log\tau.$$
Thus $$u^2\sim 2\log (T_1 T_2)$$
and $$\log u=\frac{1}{2}\log 2 +\frac{1}{2}\log\log (T_1T_2) +o(1).$$
Moreover
\begin{equation}\label{22}
u^2 = 2\log (T_1T_2)
+ \left(\frac{2}{\alpha_1}+\frac{2}{\alpha_2}-1\right)\log\log (T_1T_2)
- 4\log\tau
+ 2\log\left(\frac{\mathcal{H}_{\alpha_1}\mathcal{H}_{\alpha_2}}{2\sqrt{\pi}} 2^{1/\alpha_1+1/\alpha_2}\right)+o(1).
\end{equation}

For $T>0$ put $\delta_T=\sup_{\varepsilon\leqslant \max\cN{(|s|,|t|)} \leqslant T} \max(|r(s,t)|,\varrho_T(s,t))$.
It is straightforward to see that there exists $\delta<1$ such that for sufficiently large $T$ we get
$$\delta_T=\sup_{\varepsilon\leqslant \cN{\max(|s|,|t|)} \leqslant T} \max(|r(s,t)|,\varrho_T(s,t))<\delta<1,$$
since $\delta_T$ is decreasing in $T$ for large $T$.
Let $\beta$ be such that $0<\beta<\frac{1-\delta}{1+\delta}$.
Divide $\cN{Q:=[-T_1,T_1]\times[-T_2,T_2]-(-\varepsilon,\varepsilon)^2}$ into two subsets:
\begin{eqnarray*}
S^*\!&:=& \{(s,t)\in Q: \cN{|s|}\leqslant T_1^\beta, \cN{|t|}\leqslant T_2^\beta\},\\
S\; &:=& Q \; - \; S^*.
\end{eqnarray*}
Firstly, we show that
\begin{equation}\label{AAA}
\frac{T_1T_2}{q_1q_2}\sum_{(jq_1,kq_2)\in S^*} \rho_{T_{\max}}(jq,kq)\exp\left(-\frac{u^2}{1+\max(|r(jq,kq)|,\varrho_{T_{\max}}(jq,kq))}\right)\to 0,
\end{equation}
as $u\to\infty$.
By (\ref{22}) there exists a constant $K$ such that
$\exp(-u^2/2)\leqslant \frac{K}{T_1T_2}$. Applying the fact that $u^2\sim 2\log(T_1T_2)$ and $u^{2/\alpha_1}q_1=u^{2/\alpha_2}q_2=a$, for $u$ large enough, we obtain
\begin{eqnarray*}
\lefteqn{\frac{T_1T_2}{q_1q_2}\sum_{(jq_1,kq_2)\in S^*} \rho_{T_{\max}}(jq,kq)\exp\left(-\frac{u^2}{1+\max(|r(jq,kq)|,\varrho_{T_{\max}}(jq,kq))}\right)}\\
&&\leqslant \frac{T_1T_2}{q_1q_2}\left(\frac{\cN{2}T_1^\beta}{q_1}+1\right)\left(\frac{\cN{2}T_2^\beta}{q_2}+1\right)\exp\left(-\frac{u^2}{1+\delta}\right)
\sim \cN{4}\frac{(T_1T_2)^{\beta+1}}{q_1^2q_2^2}\left(\exp\left(-\frac{u^2}{2}\right)\right)^{\frac{2}{1+\delta}}\\
&& \leqslant \cN{4}K^{\frac{2}{1+\delta}}\frac{(T_1T_2)^{\beta+1-\frac{2}{1+\delta}}}{q_1^2q_2^2}
\sim \frac{2^{2/\alpha_1+2/\alpha_2\cN{+2}}K^{\frac{2}{1+\delta}}}{a^4} \big(\log(T_1T_2)\big)^{2/\alpha_1+2/\alpha_2}(T_1T_2)^{\beta-\frac{1-\delta}{1+\delta}}.
\end{eqnarray*}
Since we choose $\beta<\frac{1-\delta}{1+\delta}$, then (\ref{AAA}) holds.

To complete the proof it suffices to show that, as $u\to\infty$,
\begin{equation}\label{BBB}
\frac{T_1T_2}{q_1q_2}\sum_{(jq_1,kq_2)\in S} \rho_{T_{\max}}(jq_1,kq_2)\exp\left(-\frac{u^2}{1+\max(|r(jq_1,kq_2)|,\varrho_{T_{\max}}(jq_1,kq_2))}\right)\to 0.
\end{equation}
In order to do it observe that there exist constants $C>0$ and $K>0$ such that
$$\max\big(|r(s,t)|, \varrho_{T_{\max}}(s,t)\big)\cdot\log\left(\sqrt{s^2+t^2}\right)\leqslant K$$
for all $u$ sufficiently large and $(s,t)$ satisfying $C\leqslant \cN{\max(|s|,|t|)}\leqslant T_{\max}$. Put $T_{\min}:=\min(T_1,T_2)$.
Since $T_{\min}^\beta>C$ for $u$ large enough, then for $(jq_1,kq_2)$
such that $\cN{\max(|jq_1|,|kq_2|)}\geqslant T_{\min}^\beta$ we have
$$\max\big(|r(jq_1,kq_2)|, \varrho_{T_{\max}}(\cN{jq_1,kq_2})\big)\leqslant \frac{K}{\log T_{\min}^\beta}.$$
Hence
$$\exp\!\left(-\frac{u^2}{1+\max\left(|r(jq_1,kq_2)|,\varrho_{T_{\max}}(jq_1,kq_2)\right)}\right)
\!\! \leqslant \!
\exp\!\left(-\frac{u^2}{1+\frac{K}{\log T_{\min}^\beta}}\right)
\!\!\leqslant\!
\exp\!\left(-u^2\left(1-\frac{K}{\log T_{\min}^\beta}\right)\right),$$
which implies the following chain of inequalities
\begin{eqnarray*}
\lefteqn{\frac{T_1T_2}{q_1q_2}\sum_{(jq_1,kq_2)\in S} \rho_{T_{\max}}(jq_1,kq_2)
\exp\left(-\frac{u^2}{1+\max\left(|r(jq_1,kq_2)|,\varrho_{T_{\max}}(jq_1,kq_2)\right)}\right)}\\
&&\leqslant\frac{T_1T_2}{q_1q_2}\sum_{(jq_1,kq_2)\in S} \left|r(jq_1,kq_2)-\frac{r}{\log T_{\max}}\right| \exp\left(-u^2\left(1-\frac{K}{\log T_{\min}^\beta}\right)\right)\\
&&\leqslant \cN{4}\frac{T_1^2T_2^2}{q_1^2q_2^2}\exp\left(-u^2\left(1-\frac{K}{\log T_{\min}^\beta}\right)\right)\frac{1}{\log T_{\min}^{\beta}}\times
\frac{q_1q_2\log T_{\min}^{\beta}}{T_1T_2}\sum_{(jq_1,kq_2)\in S}\left|r(jq_1,kq_2)-\frac{r}{\log T_{\max}}\right|\\
&&=:I_1\times I_2.
\end{eqnarray*}

Firstly, we show that factor $I_1$ is bounded. Indeed, using that
$$u^2=2\log(T_1T_2)+\left(\frac{2}{\alpha_1}+\frac{2}{\alpha_2}-1\right)\log\log (T_1T_2)+O(1),$$
there exists a constant $K'$ such that for $u$ large enough
$$-u^2\left(1-\frac{K}{\log T_{\min}^\beta}\right) =-u^2+K\frac{2\log(T_1T_2)+\left(\frac{2}{\alpha_1}+\frac{2}{\alpha_2}-1\right)\log\log (T_1T_2) +O(1)}{\log T_{\min}^\beta}
\leqslant-u^2+K'.$$
The last inequality follows from the fact that $\frac{\log(T_1T_2)}{\log T_{\min}^\beta} \to 2/\beta$.
Moreover,
$$\exp\left(-u^2\left(1-\frac{K}{\log T_{\min}^\beta}\right)\right)\leqslant K''\exp(-u^2)\leqslant K'''(T_1T_2)^{-2}(\log (T_1T_2))^{1-2/\alpha_1-2/\alpha_2},$$
for some constants $K''$, $K'''$.
Using that $u^2\sim 2\log (T_1T_2)$ and $u^{2/\alpha_1}q_1=u^{2/\alpha_2}q_2=a$, we conclude that
\begin{eqnarray*}
I_1 &\leqslant& \cN{4}\frac{T_1^2T_2^2}{q_1^2q_2^2}\exp\left(-u^2\left(1-\frac{K}{\log T_{\min}^\beta}\right)\right) \frac{1}{\log T_{\min}^{\beta}}\\
&\leqslant& \cN{4}\frac{T_1^2T_2^2}{q_1^2q_2^2} K'''(T_1T_2)^{-2}(\log (T_1T_2))^{1-2/\alpha_1-2/\alpha_2}\frac{1}{\log T_{\min}^\beta}\\
&=& \cN{4}K'''2^{2/\alpha_1+2/\alpha_2}\frac{1}{a^4}
(\log (T_1T_2))^{2/\alpha_1+2/\alpha_2}(\log(T_1T_2))^{1-2/\alpha_1-2/\alpha_2} \frac{1}{\log T_{\min}^\beta}
\sim \frac{K'''2^{2/\alpha_1+2/\alpha_2\cN{+3}}}{a^4 \beta},
\end{eqnarray*}
which proves that $I_1$ is bounded.

In the next step we show that $I_2$ tends to $0$ as $u\to\infty$.
Observe that
\begin{eqnarray*}
I_2&=&\frac{q_1q_2\log T_{\min}^{\beta}}{T_1T_2}\sum_{(jq_1,kq_2)\in S}\left|r(jq_1,kq_2)-\frac{r}{\log T_{\max}}\right|\\
& \leqslant & \frac{q_1q_2}{T_1T_2}\sum_{(jq_1,kq_2)\in S}\left|r(jq_1,kq_2)\log(\sqrt{(jq_1)^2+(kq_2)^2}-r\right|\\
& + & \beta r \frac{q_1q_2}{T_1T_2} \sum_{(jq_1,kq_2)\in S}\left|1-\frac{\log T_{\max}}{\log(\sqrt{(jq_1)^2+(kq_2)^2}}\right|=:J_1+J_2.
\end{eqnarray*}
Combining \textbf{A3} with the fact that $a_n\to a$ implies the convergence $(a_1+a_2+\ldots+a_n)/n\to a$, as $n\to\infty$
(see \cite{RUDIN}),
we conclude that $J_1$ tends to $0$, as $u\to\infty$.
Additionally, see \cite[p. 135]{LEAD},
\begin{eqnarray*}
J_2&\leqslant&\frac{\beta r}{\log T_{\min}^\beta}  \frac{q_1q_2}{T_1T_2} \sum_{(jq_1,kq_2)\in S}
\left| \log\sqrt{(jq_1)^2+(kq_2)^2} - \log T_{\max}\right|\\
& =&\frac{r}{\log T_{\min}}\frac{q_1q_2}{T_1T_2}
\sum_{(jq_1,kq_2)\in S}
\left|\log\left(\frac{\sqrt{(jq_1)^2+(kq_2)^2}}{T_{\max}}\right)\right|
\end{eqnarray*}
Suppose that $T_{\max}=T_1$. Then
\begin{eqnarray*}
\frac{q_1q_2}{T_1T_2}
\sum_{(jq_1,kq_2)\in S}
\left|\log\left(\frac{\sqrt{(jq_1)^2+(kq_2)^2}}{T_{\max}}\right)\right|
=\frac{q_1q_2}{T_1T_2}
\sum_{(jq_1,kq_2)\in S}
\left|\log\left(\sqrt{\left(\frac{jq_1}{T_1}\right)^2+\left(\frac{kq_2}{T_2}\right)^2\left(\frac{T_2}{T_1}\right)^2}\right)\right|\\
\leqslant
\frac{q_1q_2}{T_1T_2}
\sum_{(jq_1,kq_2)\in S}
\left(\left|\log\left(\sqrt{\left(\frac{jq_1}{T_1}\right)^2+\left(\frac{kq_2}{T_2}\right)^2}\right)\right|
+ \left|\log\left|\frac{jq_1}{T_1}\right|\right|\right).
\end{eqnarray*}
Hence
\begin{eqnarray*}
J_2\leqslant \frac{r}{\log T_{\min}} O\left(\int_{\cN{-1}}^1 \int_{\cN{-1}}^1 \left|\log(\sqrt{x^2+y^2})\right|dxdy
+ \int_{\cN{-1}}^1|\log |x||dx \right)
\end{eqnarray*}
and (\ref{BBB}) holds. The combination of (\ref{AAA}) with (\ref{BBB}) completes the proof.

\end{pf}

\begin{lm}\label{LEMACIK1}
Let $q_1=q_1(u)=au^{-2/{\alpha_1}}$, $q_2=q_2(u)=au^{-2/{\alpha_2}}$ and suppose that $T=T(u)\to\infty$, as $u\to\infty$.
Then, providing that conditions \textbf{A1} and \textbf{A2} are fulfilled, there exists $\varepsilon>0$ such that
\begin{eqnarray*}
\lefteqn{\frac{m(u)}{q_1q_2}\sum_{0<\max(\cN{|jq_1|,|kq_2|}) <\varepsilon }
\Bigg[(1-r(jq_1,kq_2))\frac{r}{\log T} \left(1-\left(r(jq_1,kq_2)+(1-r(jq_1,kq_2))
\frac{r}{\log T}\right)^2\right)^{-1/2} }\\ &&\quad\quad\quad\quad\quad\quad\quad\quad\quad\quad\quad\quad\quad\quad\quad\quad
\times\exp\left(-\frac{u^2}{1+r(jq_1,kq_2)+(1-r(jq_1,kq_2)) \frac{r}{\log T}}\right)\Bigg]\to 0,
\end{eqnarray*}
as $u\to\infty$.
\end{lm}
\begin{pf}
Firstly, note that for $\varepsilon>0$ small enough
\begin{equation}\label{APROKSYMACJA}
\frac{1}{2}(\cN{|s|}^{\alpha_1} + \cN{|t|}^{\alpha_2})\leqslant 1-r(s,t)\leqslant 2(\cN{|s|}^{\alpha_1} + \cN{|t|}^{\alpha_2}),
\end{equation}
for $0\leqslant \cN{\max(|s|,|t|)}<\varepsilon$, due to \textbf{A1}. Thus for $u$ large, $\varepsilon$ small enough and $0<\max \cN{(|jq_1|,|kq_2|)}<\varepsilon$ we have
\begin{eqnarray*}
\lefteqn{\left(1-\left(r(jq_1,kq_2)+(1-r(jq_1,kq_2))
\frac{r}{\log T}\right)^2\right)^{-1/2}}\\
&\leqslant&
\left(1-\left(r(jq_1,kq_2)+(1-r(jq_1,kq_2))
\frac{r}{\log T}\right)\right)^{-1/2}
=\left((1-r(jq_1,kq_2))
\left(1-\frac{r}{\log T}\right)\right)^{-1/2}\\
&\leqslant&
\left(\frac{|jq_1|^{\alpha_1}+|kq_2|^{\alpha_2}}{4}\right)^{-1/2}
\leqslant \left(\frac{\max\left(|jq_1|^{\alpha_1},|kq_2|^{\alpha_2}\right)}{4}\right)^{-1/2}\leqslant \left(\frac{\min\left(q_1^{\alpha_1},q_2^{\alpha_2}\right)}{4}\right)^{-1/2}=Ku,
\end{eqnarray*}
for some constant $K>0$.
Combining the above inequality with (\ref{APROKSYMACJA}) and definitions of $m(u)$, $q_1$ and $q_2$ we obtain

\begin{eqnarray*}
\lefteqn{\frac{m(u)}{q_1q_2}\sum_{0<\max\cN{(|jq_1|,|kq_2|)} <\varepsilon }
\Bigg[(1-r(jq_1,kq_2))\frac{r}{\log T} \left(1-\left(r(jq_1,kq_2)+(1-r(jq_1,kq_2))
\frac{r}{\log T}\right)^2\right)^{-1/2}}\\
&&\quad\quad\quad\quad\quad\quad\quad\quad\quad\quad\quad\quad\quad\quad\quad\quad\quad \times \exp\left(-\frac{u^2}{1+r(jq_1,kq_2)+(1-r(jq_1,kq_2)) \frac{r}{\log T}}\right)\Bigg]\\
&\leqslant& K'ue^{u^2/2}\sum_{0<\max\cN{(|jq_1|,|kq_2|)} <\varepsilon }
\Bigg[\left(|jq_1|^{\alpha_1}+|kq_2|^{\alpha_2}\right)(1+\delta)\frac{ru}{\log T} \\
&&\quad\quad\quad\quad\quad\quad\quad\quad\quad\quad\quad\quad\quad\quad\quad\quad\quad  \times \exp\left(-\frac{u^2}{2-\left(|jq_1|^{\alpha_1}+|kq_2|^{\alpha_2}\right) (1-\delta - \frac{r(1+\delta)}{\log T})}\right)\Bigg]\\
&=& K'\frac{ru^2}{\log T}\sum_{0<\max\cN{(|jq_1|,|kq_2|)} <\varepsilon }
\Bigg[\left(|jq_1|^{\alpha_1}+|kq_2|^{\alpha_2}\right)(1+\delta) \\
&&\quad\quad\quad\quad\quad\quad\quad\quad\quad\quad\quad\quad\quad\quad\quad\quad\quad
\times\exp\left(-\frac{u^2 \left(|jq_1|^{\alpha_1}+|kq_2|^{\alpha_2}\right) (1 - \delta -\frac{r(1+\delta)}{\log T})}{4-2\left(|jq_1|^{\alpha_1}+|kq_2|^{\alpha_2}\right) (1 - \delta - \frac{r(1+\delta)}{\log T})}\right)\Bigg]\\
&\leqslant& K'\frac{ru^2}{\log T}(1+\delta)\frac{8}{u^2}\sum_{0<\max\cN{(|jq_1|,|kq_2|)} <\varepsilon}
\frac{u^2}{8}\left(|jq_1|^{\alpha_1}+|kq_2|^{\alpha_2}\right)
\exp\left(-\frac{u^2 \left(|jq_1|^{\alpha_1}+|kq_2|^{\alpha_2}\right)}{8}\right)\\
&=& \frac{8rK'(1+\delta)}{\log T}\sum_{0<\max\cN{(|jq_1|,|kq_2|)} <\varepsilon }
\left(\frac{|aj|^{\alpha_1}}{8}+\frac{|ak|^{\alpha_2}}{8}\right)
\exp\left(-\left(\frac{|aj|^{\alpha_1}}{8}+\frac{|ak|^{\alpha_2}}{8}\right)\right)\\
&=& O\left(\frac{K''}{\log T}\int_{\cN{-\infty}}^{\infty} \int_{\cN{-\infty}}^{\infty} \left(|x|^{\alpha_1}+|y|^{\alpha_2}\right)e^{-(|x|^{\alpha_1}+|y|^{\alpha_2})}\;dxdy\right),
\end{eqnarray*}
as $u\to\infty$.
Since $\log T(u)\to\infty$, as $u\to\infty$, and an integral in the last statement is finite, the proof is completed.
\end{pf}

\subsection{Proof of Theorem \ref{LEM2WYM}}\label{s.proof.th2}
{\bf Proof of (i)}.
Let $\{X^{(j,k)}(s,t)\}_{j,k}$ be independent copies of $X(s,t)$ and let $\eta(s,t)$ be such that $\eta(s,t)=X^{(j,k)}(s,t)$ for $(s,t)\in [j-1,j)\times[k-1,k)$.
For a fixed $T$ we define a Gaussian random field $Y_T$ as follows
\begin{equation}\label{DEF_Y_T}
Y_T(s,t):=\left(1-\frac{r}{\log T}\right)^{1/2}\eta(s,t)+\left(\frac{r}{\log T}\right)^{1/2}\mathcal{W}, \text{ \quad for $(s,t)\in[0,T]^2$},
\end{equation}
where $\mathcal{W}$ is \pE{an $N(0,1)$ random variable} independent of $\eta(s,t)$. Then the covariance of $Y_T$ equals
\[ Cov(Y_T(s_0,t_0),Y_T(s_0+s,t_0+t)) = \left\{ \begin{array}{ll}
r(s,t)+(1-r(s,t))\frac{r}{\log T},
&\text{ when }  [s_0]=[s_0+s], [t_0]=[t_0+t];\\
\frac{r}{\log T},
& \text{ otherwise }, \end{array} \right.
\]
for all $s_0,t_0,s,t\geqslant 0$.

Let $n_x:=\left\lfloor x m_1(u) \right\rfloor$ and $n_y:=\left\lfloor y m_2(u)\right\rfloor$. Since
\begin{eqnarray*}\label{OBS}
\lefteqn{P\left(\sup_{(s,t)\in[0,n_x+1]\times[0,n_y+1]}X(s,t)\leqslant u\right)}\\
&& \leqslant P\left(\sup_{(s,t)\in \left[0,xm_1(u)\right]\times\left[0,ym_2(u)\right]}X(s,t)\leqslant u\right) \leqslant P\left(\sup_{(s,t)\in[0,n_x]\times[0,n_y]}X(s,t)\leqslant u\right),
\end{eqnarray*}
we focus on the asymptotics of $P\left(\sup_{(s,t)\in[0,n_x]\times[0,n_y]}X(s,t)\leqslant u\right)$, as $u\to\infty$.
Let $\varepsilon>0$. Divide $[0,n_x]\times[0,n_y]$ into $n_xn_y$ unit squares and then split them into subsets $I_{l,m}^*$ and $I_{l,m}$ as follows
\begin{eqnarray*}
I_{l,m}  &=&[(l-1)+\varepsilon,l]\times[(m-1)+\varepsilon,m],\\
I_{l,m}^*&=&[l-1,l]\times[m-1,m] \;-\; I_{l,m},
\end{eqnarray*}
where $l=1,\ldots,n_x$, $m=1,\ldots,n_y$.

\textbf{Step 1.}
In the first step we prove that
\begin{equation}\label{NIER1}
\lim_{u\to\infty} \left|P\left(\sup_{(s,t)\in[0,n_x]\times[0,n_y]}X(s,t)\leqslant u\right)-P\left(\sup_{(s,t)\in\bigcup_{l=1}^{n_x}\bigcup_{m=1}^{n_y}I_{l,m}}X(s,t)\leqslant u\right)\right|\leqslant\rho_1(\varepsilon),
\end{equation}
uniformly for $(x,y)\in[A_0,A_\infty]^2$ with $\rho_1(\varepsilon)\to 0$ as $\varepsilon\to 0$.
This is a consequence of the following sequence of inequalities
\begin{eqnarray*}
\lefteqn{0\leqslant P\left(\sup_{(s,t)\in\bigcup_{l=1}^{n_x}\bigcup_{m=1}^{n_y}I_{l,m}}X(s,t)\leqslant u\right)-P\left(\sup_{(s,t)\in[0,n_x]\times[0,n_y]}X(s,t)\leqslant u\right)} \\
&& \leqslant n_xn_y P\left(\sup_{(s,t)\in I_{1,1}^*}X(s,t)>u\right)\leqslant A_\infty^2 m(u)P\left(\sup_{(s,t)\in I_{1,1}^*}X(s,t)>u\right)=(2\varepsilon-\varepsilon^2) A_\infty^2(1+o(1)),
\end{eqnarray*}
as $u\to\infty$,
since
$$P\left(\sup_{(s,t)\in I_{1,1}^*}X(s,t)>u\right)=\frac{2\varepsilon-\varepsilon^2}{m(u)}(1+o(1)),$$
as $u\to\infty$,
by \cite[\tr{Theorem} 7.1 ]{PITERBARG}.

\textbf{Step 2.}
Let $a>0$ and $q_1=q_1(u):=au^{-\alpha_1/2}$, $q_2=q_2(u):=au^{-\alpha_2/2}$. We show that
\begin{eqnarray}\label{NIER2}
\lim_{u\to\infty} \left|P\left(X(s,t)\leqslant u; (s,t)\in \bigcup_{l=1}^{n_x}\bigcup_{m=1}^{n_y}I_{l,m}\right)
- P\left(X(jq_1,kq_2)\leqslant u; (jq_1,kq_2)\in \bigcup_{l=1}^{n_x}\bigcup_{m=1}^{n_y}I_{l,m}\right) \right| \nonumber\\
\leqslant \rho_2(a),
\end{eqnarray}
uniformly for $(x,y)\in[A_0,A_\infty]^2$, with $\rho_2(a)\to 0$ as $a\to 0$.
Indeed, (\ref{NIER2}) follows from the fact that
\begin{eqnarray}\label{ARG}
0&\leqslant& P\left(X(s,t)\leqslant u; (s,t)\in \bigcup_{l=1}^{n_x}\bigcup_{m=1}^{n_y}I_{l,m}\right)
-P\left(X(jq_1,kq_2)\leqslant u; (jq_1,kq_2)\in \bigcup_{l=1}^{n_x}\bigcup_{m=1}^{n_y}I_{l,m}\right)\nonumber\\
&& \leqslant n_xn_y \max_{l,m}\left[P\left(X(jq_1,kq_2)\leqslant u; (jq_1,kq_2)\in I_{l,m}\right)-P\left(\sup_{(s,t)\in I_{l,m}} X(s,t)\leqslant u)\right)\right]\nonumber\\
&& \leqslant n_xn_y(1-\varepsilon)^2\left(\frac{\rho(a)}{m(u)}+o\left(\frac{1}{m(u)}\right)\right)\\
&& \leqslant A_\infty^2 \rho(a) + A_\infty^2m(u)o\left(\frac{1}{m(u)}\right)\to A_\infty^2 \rho(a),\nonumber
\end{eqnarray}
as $u\to\infty$ with $\rho(a)\to 0$ as $a\to 0$. Inequality (\ref{ARG}) is due to Lemma \ref{LEMAT1}.

\textbf{Step 3.}
In this step we show that for $T=T(u):=\max(A_{\infty}m_1(u),A_{\infty}m_2(u))$ we have
\begin{equation}\label{NIE3}
\left|P\left(X(jq_1,kq_2)\leqslant u; (jq_1,kq_2)\in \bigcup_{l=1}^{n_x}\bigcup_{m=1}^{n_y}I_{l,m}\right)
\!\!- \!\!P(Y_{T}(jq_1,kq_2)\leqslant u; (jq_1,kq_2)\in \bigcup_{l=1}^{n_x}\bigcup_{m=1}^{n_y}I_{l,m})\right|\to 0,
\end{equation}
as $u\to\infty$, uniformly for $(x,y)\in[A_0,A_\infty]^2$.

Indeed, note that for sufficiently large $T$ we have
\begin{eqnarray*}
\big|Cov(X(jq_1,kq_2),X(j'q_1,k'q_2))-Cov(Y_T(jq_1,kq_2),Y_T(j'q_1,k'q_2))\big|
&\leqslant& \rho_T((j-j')q_1, (k-k')q_2),\\
\big|Cov(Y_T(jq_1,kq_2),Y_T(j'q_1,k'q_2))\big|
&\leqslant& \varrho_T((j-j')q_1, (k-k')q_2),
\end{eqnarray*}
for functions $\rho_T$ and $\varrho_T$ defined by (\ref{FUN}).

Moreover, for \cN{small $\varepsilon>0$} and $(jq_1,kq_2),(j'q_1,k'q_2)\in \bigcup_{l=1}^{n_x}\bigcup_{m=1}^{n_y}I_{l,m}$
satisfying $\max(|j-j'|q_1,|k-k'|q_2)<\varepsilon$ we get
\begin{equation*}
\big|Cov(X(jq_1,kq_2),X(j'q_1,k'q_2))-Cov(Y_T(jq_1,kq_2),Y_T(j'q_1,k'q_2))\big|= (1-r((j-j')q_1,(k-k')q_2))\frac{r}{\log T}
\end{equation*}
and
\begin{eqnarray*}
\lefteqn{\max\left(|Cov(X(jq_1,kq_2),X(j'q_1,k'q_2))|,
|Cov(Y_T(jq_1,kq_2),Y_T(j'q_1,k'q_2)|\right
)}\\
&=&Cov(Y_T(jq_1,kq_2),Y_T(j'q_1,k'q_2)) \\
&=&r((j-j')q_1,(k-k')q_2)+(1-r((j-j')q_1,(k-k')q_2)) \frac{r}{\log T}.
\end{eqnarray*}
Let $\delta_T=\sup\{\max(|r(s,t)|,\varrho_T(s,t));\;\max\cN{(|s|,|t|)}\geqslant \varepsilon\}$. Observe that $\delta_T<\delta<1$ for sufficiently large $T$.
Applying \cite[Theorem 4.2.1]{LEAD} we get
\begin{eqnarray*}
\lefteqn{\left|P\left(X(jq_1,kq_2)\leqslant u; (jq_1,kq_2)\in \bigcup_{l=1}^{n_x}\bigcup_{m=1}^{n_y}I_{l,m}\right)
- P\left(Y_T(jq_1,kq_2)\leqslant u; (jq_1,kq_2)\in \bigcup_{l,m}I_{l,m}\right)\right|}\\
&\leqslant &
\frac{1}{4\pi}\frac{n_xn_y}{q_1q_2}\sum_{0<\max\cN{(|jq_1|,|kq_2|)} <\varepsilon }
\Bigg[(1-r(jq_1,kq_2))\frac{r}{\log T} \\
&& \times \left(1-\left(r(jq_1,kq_2)+(1-r(jq_1,kq_2))
\frac{r}{\log T}\right)^2\right)^{-1/2}
\exp\left(-\frac{u^2}{1+r(jq_1,kq_2)+(1-r(jq_1,kq_2))
\frac{r}{\log T}}\right)\Bigg]\\
&+&
\frac{1}{4\pi}(1-\delta^2)^{-1/2}\frac{n_xn_y}{q_1q_2}\sum_{(jq_1,kq_2) \in \cN{[-n_x,n_x]\times[-n_y,n_y]-(-\varepsilon,\varepsilon)^2}}
\Bigg[ \rho_T(jq_1,kq_2) \\
&&\times\exp\left(-\frac{u^2}{1+\max(|r(jq_1,kq_2)|,\varrho_T(jq_1,kq_2))}\right)\Bigg]
\\
&\leqslant&
\frac{1}{4\pi}\frac{A_\infty^2m(u)}{q_1q_2}\sum_{0<\max\cN{(|jq_1|,|kq_2|)} <\varepsilon }
\Bigg[(1-r(jq_1,kq_2))\frac{r}{\log T}\\
&& \times \left(1-\left(r(jq_1,kq_2)+(1-r(jq_1,kq_2))
\frac{r}{\log T}\right)^2\right)^{-1/2}
\exp\left(-\frac{u^2}{1+r(jq_1,kq_2)+(1-r(jq_1,kq_2)) \frac{r}{\log T}}\right)\Bigg]\\
&+&
\frac{1}{4\pi} (1-\delta^2)^{-1/2}\frac{A_\infty^2m(u)}{q_1q_2}\times
\sum_{(jq_1,kq_2)\in \cN{[-A_\infty m_1(u),A_\infty m_1(u)]\times [-A_\infty m_2(u),A_\infty m_2(u)] - (-\varepsilon,\varepsilon)^2}}
\Bigg[\rho_T(jq_1,kq_2)\\
&&\times
\exp\left(-\frac{u^2}{1+\max(|r(jq_1,kq_2)|,\varrho_T(jq_1,kq_2))}\right)\Bigg]\\
&=:& I_1 + I_2.
\end{eqnarray*}
Observe that, due to Lemma \ref{LEMACIK1}, $I_1$ tends to 0 as $u\to\infty$. Analogously, by Lemma \ref{LEMACIK2}, $I_2$ tends to $0$
as $u\to\infty$. Hence we have shown (\ref{NIE3}).

\textbf{Step 4.}
By definition of the random field $Y_T$, we have
\begin{eqnarray}\label{AAAA}
\lefteqn{P\left(Y_T(jq_1,kq_2)\leqslant u; (jq_1,kq_2)\in \bigcup_{l,m}I_{l,m}\right)} \nonumber\\
&=& P\left(\left(1-\frac{r}{\log T}\right)^{1/2}\eta(jq_1,kq_2)+\left(\frac{r}{\log T}\right)^{1/2}\mathcal{W}\leqslant u;\;(jq_1,kq_2)\in \bigcup_{l,m}I_{l,m}\right) \nonumber \\
&=& P\left( \left(1-\frac{r}{\log T}\right)^{1/2}
\sup_{(jq_1,kq_2)\in \bigcup_{l,m}I_{l,m}} \eta(jq_1,kq_2) + \left(\frac{r}{\log T}\right)^{1/2}\mathcal{W} \leqslant u \right)\nonumber \\
&=&\int_{-\infty}^{\infty}P\left(\sup_{(jq_1,kq_2)\in\bigcup_{l,m}I_{l,m}}\eta(jq_1,kq_2)\leqslant \frac{u-(r/\log T)^{1/2}z}{(1-r/\log T)^{1/2}}\right)d\Phi(z).
\end{eqnarray}
Then for any $z\in \mathbb{R}$
\begin{eqnarray*}
u_z&:=&\frac{u-(r/\log T)^{1/2}z}{(1-r/\log T)^{1/2}}\\
&=&\left(u-(r/\log T)^{1/2}z\right)\left( 1+\frac{1}{2}(r/\log T)+o(r/\log T)\right) \\
&=&u+\frac{-2\sqrt{r}z+2r}{u}+o(1/u),
\end{eqnarray*}
as $u\to\infty$,
and thus
\tb{
\begin{eqnarray*}
\frac{1}{m(u_z)}&=&
\frac{\exp(-2r+2\sqrt{r}z)}{m(u)}(1+o(1)).
\end{eqnarray*}
}
Hence, we get
\tb{
\begin{eqnarray}\label{BBBB}
P\left(\sup_{\cN{(jq_1,kq_2)}\in\bigcup_{l,m}I_{l,m}}\eta(jq_1,kq_2)\leqslant u_z\right)&=&
\prod_{l,m}P\left(\sup_{(jq_1,kq_2)\in I_{l,m}}X(jq_1,kq_2)\leqslant u_z\right)\nonumber\\
&=& P\left(\sup_{(s,t)\in [0,1]^2}X(s,t)\leqslant u_z\right)^{n_xn_y}(1+o(1))\nonumber\\
&=& \left(1-\frac{1}{m(u_z)}\right)^{xym(u)}(1+o(1))\nonumber\\
&=&\exp(-xy\exp(-2r+2\sqrt{r}z))(1+o(1)),
\end{eqnarray}
}
as $u\to\infty$, uniformly for $(x,y)\in [A_0,A_{\infty}]^2$.
Combining (\ref{NIER1}), (\ref{NIER2}), (\ref{NIE3}), (\ref{AAAA}) and (\ref{BBBB})
and passing with $\varepsilon\to 0$ and $a\to 0$,
we conclude that the proof of (i) is completed.
\\
\\
\tb{
{\bf Proof of (ii)}. Let $\mathcal{T}\subset\rr^2$ be Jordan-measurable with Lebesgue measure ${\rm mes}(\mathcal{T})>0$.
For given $\varepsilon>0$,
let
$\mathcal{L}_\varepsilon,\mathcal{U}_\varepsilon\subset \rr^2$ be {\it simple sets} (i.e. finite sums of disjoint
rectangles of the form $[a_1,b_1)\times[a_2,b_2)$)
such that
 $\mathcal{L}_\varepsilon\subset \mathcal{T} \subset \mathcal{U}_\varepsilon$
and ${\rm mes}(\mathcal{L}_\varepsilon)>{\rm mes}(\mathcal{T})-\varepsilon$,
${\rm mes}(\mathcal{U}_\varepsilon)<{\rm mes}(\mathcal{T})+\varepsilon$.
Then, following line-by-line the same argument as given in the proof of
part (i) of Theorem \ref{LEM2WYM}, for
$\mathcal{T}_u=\{(x,y):(x/ m_1(u),y /m_2(u))\in\mathcal{T}\},
\mathcal{L}_{\varepsilon,u}=\{(x,y):(x/m_1(u),y/m_2(u))\in\mathcal{L_\varepsilon}\},
\mathcal{U}_{\varepsilon,u}=\{(x,y):(x/m_1(u),y/m_2(u))\in\mathcal{U_\varepsilon}\}$
we have
\[
P\left(\sup_{(s,t)\in \mathcal{L}_{\varepsilon,u}}X(s,t)\leqslant u\right)
\to
E\big(\exp(- {\rm mes}(\mathcal{L}_{\varepsilon})\exp(-2r+2\sqrt{r}\mathcal{W}))\big)
\]
and
\[
P\left(\sup_{(s,t)\in \mathcal{U}_{\varepsilon,u}}X(s,t)\leqslant u\right)
\to
E\big(\exp(- {\rm mes}(\mathcal{U}_{\varepsilon})\exp(-2r+2\sqrt{r}\mathcal{W}))\big),
\]
as $u\to\infty$.
Thus,
\[
P\left(\sup_{(s,t)\in \mathcal{T}_u}X(s,t)\leqslant u\right)
\to
E\big(\exp(- {\rm mes}(\mathcal{T})\exp(-2r+2\sqrt{r}\mathcal{W}))\big),
\]
as $u\to\infty$.
}

\subsection{Proof of Proposition \ref{p.random}}\label{sub3}

Since the proof of Proposition \ref{p.random} is analogous to proofs of Theorems 3.1-3.3
in \cite{AREN}, see also
Theorem A in \cite{TAN}, we focus only on arguments for (ii).

Let $0<A_0<A_\infty$. We have
\begin{eqnarray*}
\lefteqn{
P\left(\sup_{(s,t)\in \mathcal{B}(0,T)}X(s,t)>u\right)
=}\\
&=&
\int_0^{A_0\sqrt{m(u)}} P\left(\sup_{(s,t)\in \mathcal{B}(0,x)}X(s,t)>u\right)dF_T(x)
+
\int_{A_0\sqrt{m(u)}}^{A_\infty\sqrt{m(u)}} P\left(\sup_{(s,t)\in \mathcal{B}(0,x)}X(s,t)>u\right)dF_T(x)\\
&&+
\int_{A_\infty\sqrt{m(u)}}^\infty P\left(\sup_{(s,t)\in \mathcal{B}(0,x)}X(s,t)>u\right)dF_T(x)
=I_1+I_2+I_3.
\end{eqnarray*}
Then, for each $\varepsilon>0$, due to Remark \ref{r.2}, for sufficiently large $u$, we get
\begin{eqnarray*}
I_2
&\le&
(1+\varepsilon)
\int_{A_0}^{A_\infty } (1-
E\big(\exp(-\pi x^2\exp(\mathcal{V}_r))\big) dF_T(x\sqrt{m(u)})\\
&=&
(1+\varepsilon)
\int_{A_0}^{A_\infty}
 2\pi x
E\big(\exp(-\pi x^2\exp(\mathcal{V}_r)+\mathcal{V}_r)\big)
P(T>x\sqrt{m(u)})dx\\
&&-
(1+\varepsilon)
\left(1-
E\big(\exp(-\pi A_\infty^2\exp(\mathcal{V}_r))\big)
\right)
P(T>A_\infty \sqrt{m(u)})\\
&&+
(1+\varepsilon)
\left(1-
E\big(\exp(-\pi A_0^2\exp(\mathcal{V}_r))\big)
\right)
P(T>A_0 \sqrt{m(u)}),
\end{eqnarray*}
where
$\mathcal{V}_r=2\sqrt{r}\mathcal{W}-2r$.
Hence, using the fact that $T$ is regularly varying,
\begin{eqnarray*}
\limsup_{u\to\infty} \frac{I_2}{P(T> \sqrt{m(u)})}
&\le&
(1+\varepsilon)
2\pi
\int_{A_0}^{A_\infty}
x^{1-\lambda}E\big(\exp(-\pi x^2\exp(\mathcal{V}_r)+\mathcal{V}_r)\big)dx\\
&&-
(1+\varepsilon)
\left(1-
E\big(\exp(-\pi A_\infty^2\exp(\mathcal{V}_r))\big)
\right)
A_\infty^{-\lambda}\\
&&+
(1+\varepsilon)
\left(1-
E\big(\exp(-\pi A_0^2\exp(\mathcal{V}_r))\big)
\right)
A_0^{-\lambda}.
\end{eqnarray*}
In an analogous way we get that
\begin{eqnarray*}
\liminf_{u\to\infty} \frac{I_2}{P(T> \sqrt{m(u)})}
&\geqslant&
(1-\varepsilon)
2\pi
\int_{A_0}^{A_\infty}
x^{1-\lambda}E\big(\exp(-\pi x^2\exp(\mathcal{V}_r)+\mathcal{V}_r)\big)dx\\
&&-
(1-\varepsilon)
\left(1-
E\big(\exp(-\pi A_\infty^2\exp(\mathcal{V}_r))\big)
\right)
A_\infty^{-\lambda}\\
&&+
(1-\varepsilon)
\left(1-
E\big(\exp(-\pi A_0^2\exp(\mathcal{V}_r))\big)
\right)
A_0^{-\lambda}.
\end{eqnarray*}
Then, following the same argument as in the proof of Theorem 3.2 in \cite{AREN},
we conclude that
$I_1+I_3 =o(P(T> \sqrt{m(u)}))$
as $u\to\infty$.

Now, passing with $A_0\to0$, $A_\infty\to\infty$ and $\varepsilon\to0$,
we conclude that
\[
I_2=
2\pi
\int_{0}^{\infty}
x^{1-\lambda}E\big(\exp(-\pi x^2\exp(\mathcal{V}_r)+\mathcal{V}_r)\big)dx
P(T> \sqrt{m(u)})(1+o(1)),
\]
as $u\to\infty$.

\subsection{Proof of Proposition \ref{prop.1}}\label{SUB2}

{\bf Proof of (i)}.
\tr{Assume that {\bf A3} is satisfied with $r=0$. Then, by definition of
$\{\widetilde{X}_{j,k}\}$, it suffices to show that
for the original Gaussian field $\{X(s,t):s,t\ge0\}$}
\begin{equation}
P\left(\sup_{(s,t)\in[0,f(u)]\times[0,g(u)]}X(s,t)\leqslant z(u)\right)-P\left(\sup_{(s,t)\in[0,1]^2}X(s,t)\leqslant z(u)\right)^{f(u)g(u)}\to 0
\label{BB}
\end{equation}
as $u\to\infty$, for each function $z:\mathbb{R}_+\to\mathbb{R}$ and all pairs of functions $f,g:\mathbb{R}_+\to\mathbb{R}_+$ such that $f(u)\to\infty$ and $g(u)\to\infty$, as $u\to\infty$, and $1/C\leqslant f(u)/g(u)\leqslant C$ for some fixed $C>0$.
Observe that it suffices to consider two cases: continuous $z(u)\nearrow \infty$, as $u\to\infty$, and $z(u)<Const$.
We focus on the first case and suppose that $z(u)$ increases to infinity. Then (\ref{BB}) is equivalent to
\begin{equation}
P\left(\sup_{(s,t)\in\left[0,f^*(u)\right]\times\left[0,g^*(u)\right]}X(s,t)\leqslant u\right)
-P\left(\sup_{(s,t)\in[0,1]^2}X(s,t)\leqslant u\right)^{f^*(u)g^*(u)}\to 0,
\end{equation}
as $u\to\infty$, with $z^{-1}$ being the inverse function for $z$ and \tr{$f^*(u):=f(z^{-1}(u))$, \mbox{$g^*(u):=g(z^{-1}(u))$}.}

By (i) of Theorem \ref{LEM2WYM},
\begin{equation}\label{UN2}
P\left(\sup_{(s,t)\in\left[0,x\sqrt{m(u)}\right]\times\left[0,y\sqrt{m(u)}\right]}X(s,t)\leqslant u\right)\to e^{-xy},
\end{equation}
as $u\to\infty$, uniformly for $(x,y)\in\mathcal{F}(C):=\left\{(s,t)\in\mathbb{R}_+^2:\; 1/C\leqslant s/t \leqslant C\right\}\tr{\cup \{0,0\}},$
for an arbitrary constant $C>0$.
Moreover the uniform convergence
\begin{eqnarray}
P\left(\sup_{(s,t)\in[0,1]^2}X(s,t)\leqslant u\right)^{xy\cdot m(u)}\to e^{-xy}\label{UN1}
\end{eqnarray}
occurs on the set $\mathcal{F}(C)$.

\tr{Let $\bar{f}(u):=f\left(z^{-1}(u)\right)/\sqrt{m(u)}$
and $\bar{g}(u):=g\left(z^{-1}(u)\right)/\sqrt{m(u)}$.}
The fundamental observation is that it is sufficient to prove (\ref{BB}) for $f(u)$ and $g(u)$ satisfying the additional assumption:
$\bar{f}(u)\to a\in[0,\infty]$ and $\bar{g}(u)\to b\in[0,\infty]$, as $u\to\infty$.

Note that $1/C\leqslant f(u)/g(u)\leqslant C$ implies $1/C\leqslant \bar{f}(u)/\bar{g}(u)\leqslant C$.
Since the convergence in (\ref{UN2}) is uniform, we obtain
\begin{eqnarray*}
P\left(\sup_{(s,t)\in[0,f^*(u)]\times[0,g^*(u)]}X(s,t)\leqslant u\right)
 =P\left(\sup_{(s,t)\in[0,\bar{f}(u)\sqrt{m(u)}]\times[0,\bar{g}(u)\sqrt{m(u)}]}X(s,t)\leqslant u\right)\to e^{-ab},
\end{eqnarray*}
as $u\to\infty$. On the other hand, by (\ref{UN1}),
$$P\left(\sup_{(s,t)\in[0,1]^2}X(s,t)\leqslant u\right)^{f^*(u)g^*(u)}=P\left(\sup_{(s,t)\in[0,1]^2}X(s,t)\leqslant u\right)^{\bar{f}(u)\bar{g}(u)\cdot m(u)}\to e^{-ab},$$
as $u\to\infty$, which gives (\ref{BB}).
\\
\\
{\bf Proof of (ii).} \tr{Let us consider the case $r>0$. Note that for $\mathcal{V}_r=2\sqrt{r}\mathcal{W}-2r$ it holds that
\begin{eqnarray*}
\lefteqn{Var\left(\exp(-\exp(\mathcal{V}_r))\right)=E\left(\exp\left(-2\exp(\mathcal{V}_r)\right)\right)-
E\big(\exp(-\exp(\mathcal{V}_r))\big)^2}\\
&=& P\left( \max_{j\leqslant 2\left\lfloor\sqrt{m(u)}\right\rfloor,k\leqslant \left\lfloor\sqrt{m(u)}\right\rfloor}\widetilde{X}_{j,k}\leqslant u \right)-P\left( \max_{j,k\leqslant\left\lfloor\sqrt{m(u)}\right\rfloor}\widetilde{X}_{j,k}\leqslant u \right)^2 +o(1),
\end{eqnarray*}
due to Theorem \ref{LEM2WYM}.
By contradiction, assume that the extremal index exists and equals \tr{$\theta \in (0,1]$}. Then for any sequence $(z_n)\subset \mathbb{R}$ we have
\begin{eqnarray*}
\lefteqn{P\left(\max_{j\leqslant \left\lfloor2\sqrt{m(z_n)}\right\rfloor,k\leqslant\left\lfloor\sqrt{m(z_n)}\right\rfloor}\widetilde{X}_{j,k}\leqslant z_n \right)-P\left(\max_{j,k\leqslant\left\lfloor\sqrt{m(z_n)}\right\rfloor}\widetilde{X}_{j,k}\leqslant z_n \right)^2}\\
&=&
\left(P\left(\max_{j\leqslant 2\left\lfloor\sqrt{m(z_n)}\right\rfloor,k\leqslant\left\lfloor\sqrt{m(z_n)}\right\rfloor}\widetilde{X}_{j,k}\leqslant z_n \right)-
P\left(\widetilde{X}_{1,1}\leqslant z_n\right)^{2 m(z_n)\cdot \theta}\right)\\
&&-
\left(P\left(\max_{j,k\leqslant\left\lfloor\sqrt{m(z_n)}\right\rfloor}\widetilde{X}_{j,k}\leqslant z_n \right)^2
-
\left(P\left(\widetilde{X}_{1,1}\leqslant z_n \right)^{m(z_n)\cdot\theta}\right)^2\right)=o(1),
\end{eqnarray*}
as $n\to\infty$, which implies that $Var\left(\exp(-\exp(\mathcal{V}_r))\right)=0$.
Keeping in mind that  $ r>0$ and $\mathcal{W}$ is \pE{an $N(0,1)$ random variable},
we obtain a contradiction.}
\\
\\
{\bf Acknowledgement}: K. D\c{e}bicki was partially supported by
NCN Grant No 2011/01/B/ST1/01521 (2011-2013).
The first two authors kindly acknowledge partial support by the Swiss National
Science Foundation Grant 200021-140633/1 and by the project RARE -318984 (an FP7
a Marie Curie IRSES Fellowship)


\end{document}